\documentclass[12pt,a4paper]{amsart}
\usepackage{amssymb,amsmath,amsthm}
\usepackage{graphicx}
\usepackage{xcolor}

\usepackage{epsfig}

\long\def\onefigure#1#2{
\begin{figure*}[tbp]
\begin{center}
#1
\end{center}
\caption{#2}
\end{figure*}
}

\newcommand{\lipefig}[2]  
{\onefigure{\mbox{\psfig{file=#1.eps}}}{\label{f:#1} #2} }

\newtheorem{theorem}{Theorem}[section]
\newtheorem{lemma}{Lemma}[section]

\newtheorem{claim}{Claim}[section]

\newcommand{\al}{\alpha}
\newcommand{\be}{\beta}
\newcommand{\ga}{\gamma}

\newcommand{\bd}{\partial}
\newcommand{\eps}{\varepsilon}

\newcommand{\R}{\mathbb{R}}

\newcommand{\diam}{\textrm{diam\;}}

\newcommand{\remove}[1]{}

\numberwithin{equation}{section}

\begin{document}

\title{A special balancing game, Victor Grinberg's question}
\author{Imre B\'ar\'any}

\maketitle

\begin{center}
{\sl Dedicated to the memory of Victor Grinberg (1942--2019)}
\end{center}

\begin{abstract} In a vector balancing game, on step $k$ player 1 chooses a unit vector $v_k$ and player 2 chooses a sign $\eps_k\in \{-1,1\}$, and the position after $n$ steps is $z_n=\sum_1^n\eps_k v_k$. Player 1's target is to make $\|z_n\|$ large and player 2's target is to make it small. We consider a special case of this game.
\end{abstract}

\section{Introduction and some background}\label{sec:introd}
Balancing games were introduced by Joel Spencer \cite{Spen}. They are determined by a bounded set $V \subset \R^d$ and an origin symmetric convex body $K \subset \R^d$ that defines a norm $|\,.\,|_K$ in the usual way, that is $|z|_K=\min \{t\ge 0: z\in tK\}.$ It is played by two players called Pusher and Chooser. It proceeds through a sequence of positions $z_0=0,z_1,\ldots,z_n \in \R^d$ the following way. In step $n$ the position $z_n$ is known to both players, Pusher chooses $v_n \in V$ and tells it to Chooser who picks $\eps_n \in \{-1,1\}$. The next position is $z_{n+1}=z_n+\eps_nv_n$. Pusher's aim is to make $|z_n|_K$ large and Chooser's aim is to make it small. This is the $(V,K)$ game or the $(V,K,n)$ game if one wants to emphasize the total number of steps. We mention that in~\cite{Spen}, instead of $|z_n|_K$, the value of the game is $\max_{j\in\{1,\ldots,n\}}|z_j|_K$. As we will see later this is not a significant difference.

\medskip
A simple special case is when $V=K=B$ where $B\subset \R^2$ is the unit disk in the plane and the norm $|\,.\,|_B$ is the Euclidean one,
to be denoted by $\|\,.\,\|.$ It is easy to see that Pusher can always guarantee $\|z_n\|\ge \sqrt n$ and Chooser has a strategy giving $\|z_n\|\le \sqrt n$, see \cite{Spen}. So the value of the game $(B,B,n)$ is $\sqrt n$. Note that the same holds for the game $(S^1,B,n)$ where $S^1$ is the unit circle, the boundary of $B.$ It is easy to check that the same estimate holds (with the same strategies) when $B$ is the Euclidean ball in $\R^d$.

\medskip
In the same paper \cite{Spen} Spencer considered the game $(B^{\infty},B^{\infty})$ where $B^{\infty}$ is the unit ball of the max norm in $\R^d$ and showed that  the value of this game is between $\sqrt {1-\frac 1d}\sqrt n$ and $\sqrt {2\log 2d}\sqrt n$. Spencer asked what happens when $V \subset \R^d$ is a finite set. In this case Chooser can keep the set of positions in the set $S(V):=\{\sum_{w\in W}: W\subset V\}$, so the set of positions is bounded, and $\|z_n\|$ is not larger than the diameter of $S(V)$. A strategy for Pusher given in \cite{Bar79} guarantees that after sufficiently many rounds $z_n$ visits positions at least at distance $\frac 12 \diam S(V)$ from the origin (provided $V$ contains no collinear vectors).

\medskip
There are several related results, see the survey paper \cite{Bar}. A typical example is that given a sequence $v_1,\ldots,v_n \in C$, where $C\subset \R^d$ is the unit ball of the norm $|\,.\,|_C$, are there signs $\eps_1,\ldots,\eps_n$ such that $\max_k|\sum_1^k \eps_iv_i|_C$ is small, or in another version $|\sum_1^n \eps_iv_i|_C$ is small? It is shown in \cite{BarGr} that there exist signs such that $\max_k|\sum_1^k \eps_iv_i|_C\le 2d-1$ for every norm in $\R^d$. It is known \cite{Huafei} that if on each step in game $(C,C)$, Chooser may leave just one $\eps_i$ undecided, then she can guarantee $\max_k|\sum_1^k \eps_iv_i|_C \le 3$, here  $C$ is the unit ball of an arbitrary norm in $\subset \R^2$. (This result is implicit in \cite{BarGr}.) Swanepoel \cite{Swa} proves that if $v_1,\ldots,v_{2k+1}$ are unit vectors in some normed plane, then there are signs $\eps_1,\ldots,\eps_{2k+1}$ such that the norm of $\sum_{i=1}^{2k+1} \eps_iv_i$ is at most one. A result of Blokhuis and Chen~\cite{Blok} states that if $U=\{u_1,\ldots,u_n\}$ is a set of $n\approx \frac 12 d\log d$ (Euclidean) unit vectors in $\R^d$, then there are signs $\eps_i$ such that $\|\sum_{i=1}^{n} \eps_iv_i\|\le 1$. A striking question of Koml\'os is the following. Given  $v_1,\ldots,v_n\in \R^d$, with each $\|v_i\|\le 1$, are there signs $\eps_1,\ldots,\eps_n$ so that the maximum norm of  $\sum_1^n \eps_iv_i$ is bounded by a constant independent of $d$ and $n$? The best result in this direction is due to Banaszczyk~\cite{Ban} who proved the bound $c\sqrt{\log n}$ where $c$ is a universal constant.

\section{Main result}\label{sec:main}

The following variant of this game was proposed by Victor Grinberg, my friend and colleague and multiple coauthor, in an email \cite{Grin} in 2018. What is the value of the game $(S^*,B,n)$ where $S^*=\{(x,y)\in S^1: |y| \le |x|\}$? As $S^* \subset S^1$, Chooser's previous strategy guarantees that $\|z_n\|\le \sqrt n$. But can she achieve $\|z_n\|\le c\sqrt n$ for some $c<1$? Also, can Pusher guarantee $\|z_n\|\ge c\sqrt n$ for some $c>0$? The aim of this paper is to show that the answer to both questions is yes:

\begin{theorem}\label{th:main} In game $(S^*,B)$ Pusher has a strategy giving
\[
\|z_n\| \ge (\sqrt 2 -1)\sqrt n \mbox{ for every } n,
\]
and Chooser has a strategy showing that
\[
\|z_n\| \le  0.9722 \sqrt n +O\left(\frac {\log n}{\sqrt n}\right) \mbox{ for every } n.
\]
\end{theorem}

This means that, for large enough $n$, $\frac {\|z_n\|}{\sqrt n}$ is between 0.41 and 0.98, still a large gap between the lower and upper bounds. I think that $\lim \sup \frac {\|z_n\|}{\sqrt n}$ is close to the upper bound.

\section{Pusher's strategy in Theorem~\ref{th:main}}\label{sec:Push}

Define $B(t,r) \subset \R^2$ as the disk of radius $r>0$ centered at $(0,t)\in \R^2$. The set $K=B(1,\sqrt 2)\cap B(-1,\sqrt 2)$ is an $0$-symmetric body in $\R^2$, see Figure 1 left. We are going to show that in the game $(S^*,K)$ Pusher has a strategy giving
\begin{equation}\label{eq:push}
|z_n|_K\ge \sqrt n.
\end{equation}

\begin{figure}[h]
\centering
\includegraphics[scale=0.8]{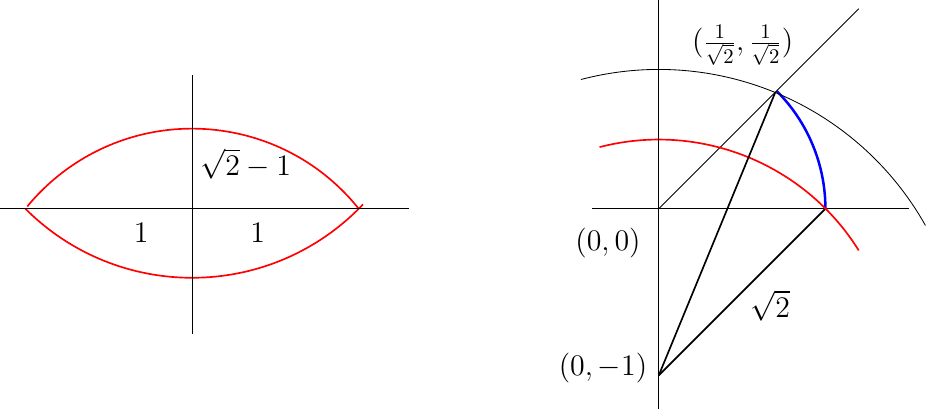}
\caption{The set $K$ and $|v_1|_K$}
\label{fig:Kv1}
\end{figure}

We show first that for every $v \in S^*$
\begin{equation}\label{eq:v_1}
1\le |v|_K\le \frac 12+ \sqrt {\frac 32}.
\end{equation}
We may assume that $v=(\cos \phi, \sin \phi)$ with $\phi \in [0,\pi/4]$.
Then $|v|_K$ is minimal when $\phi=0$, so indeed $|v|_K \ge 1,$  see Figure 1 right.
Analogously $|v|_K$ is maximal when $\phi=\pi/4$ and then $|v|_K=\frac 12+ \sqrt {\frac 32}.$

Pusher's strategy is simple: in step $n$ he computes $t_n=|z_n|_K$ and selects a unit vector $v_n\in S^*$ tangent to $t_nK$ at $z_n\in \bd (t_nK)$. Note that by the choice of $K$ such a unit vector always exists. Then $z_{n+1}=z_n+\eps_nv_n$ is exactly at distance $\sqrt{2t_n^2+1}$ from $(0,-t_n)$ or from $(0,t_n)$ see Figure 2. Thus $z_{n+1}$ is not lying in the interior of the intersection
$$
B(t_n,\sqrt {2t_n^2+1})\cap B(-t_n,\sqrt {2t_n^2+1}).
$$

\begin{figure}[h]
\centering
\includegraphics[scale=0.8]{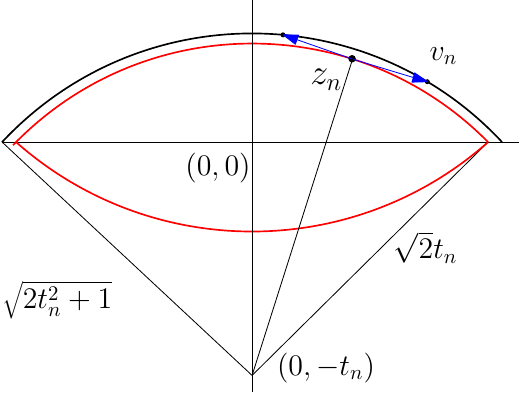}
\caption{The position $z_n$ and Chooser's choice $v_n$}
\label{fig:tn}
\end{figure}

Direct checking shows that this intersection contains $\sqrt{t_n^2+1}K$. Thus $t_{n+1}=|z_{n+1}|_K \ge \sqrt{t_n^2+1}$,
implying $t_{n+1}^2\ge t_n^2+1$. Then
\[
t_n^2\ge t_{n-1}^2+1 \ge t_{n-2}^2+2\ge \ldots \ge t_1^2+(n-1)\ge n
\]
because $t_1\ge 1$ for every choice of $v_1$, according to (\ref{eq:v_1}).

Since $B(0,\sqrt 2-1) \subset K$, inequality \ref{eq:push} shows that $\|z_n\|\ge (\sqrt 2-1)\sqrt n$ finishing the proof of the first half of Theorem~\ref{th:main}. \qed

\medskip
One can show (we omit the details) that Chooser has a strategy in game $(S^*,K)$ guaranteeing
\begin{equation}\label{eq:choo}
|z_n|_K\le  1.6002.. \sqrt n +O\left(\frac {\log n}{\sqrt n}\right).
\end{equation}

\medskip
Unfortunately this only gives $\|z_n\|\le 1.6003\sqrt n$ (for large enough $n$) much weaker than $\|z_n\|\le \sqrt n$ that we have from game $(S^1,B)$.

\medskip
The difficulty is that if $z_n=(x_n,y_n)$ satisfies $|y_n|\ge |x_n|$ for every $n$, then Pusher's strategy from game $(S^1,B)$ guarantees $\|z_n\|\ge \sqrt n$.
That is why Chooser has to force the position $z_n$ to go to the region where $|y_n|\le |x_n|$. This is explained in the next section.

\section{Chooser's strategy in game $(S^*,D(\al))$}\label{sec:Dalpha}

For simpler notation we define $D(t,r)=B(t,r)\cap B(-t,r)$, here $t,r>0$. For a fixed $\al\ge 0$ set $D(\al)=D(1,\sqrt 2(\al+1))$. Note that for $\al=0$, $D(\al)$ coincides with $K$, $D(0)=K$. We assume that from now on $\al>0$.

\begin{figure}[h]
\centering
\includegraphics[scale=0.8]{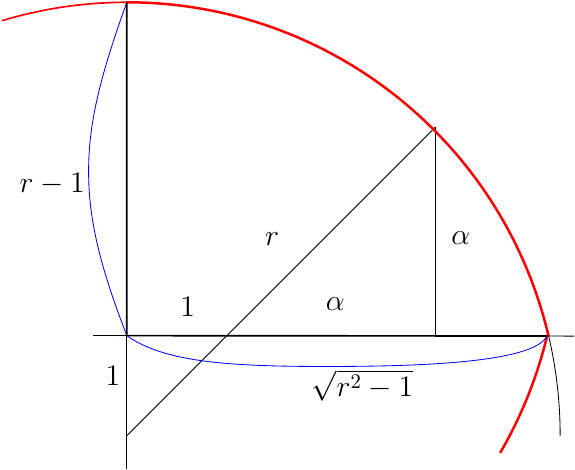}
\caption{The top right quarter of $D(\alpha)$. }
\label{fig:DB}
\end{figure}

\smallskip
$D(\al)$ is a compact convex set, symmetric with respect to the $x$ and $y$ axes, that contains the origin in its interior, see Figure~\ref{fig:DB}.
It defines a norm $|\,.\,|_{D(\al)}$.
Define further the (double) cones $C_1=\{(x,y)\in \R^2: (1+\al)|y|\le \al|x|\}$ and $C_2=\{(x,y)\in \R^2: (1+\al)|y|\ge \al|x|\}$.

\medskip
Suppose $z=(x,y)$ is the position in a step of game $(S^*,D(\al))$. We can assume (by symmetry) that $x,y\ge 0$. Set $t=t(z)\!:=|z|_{D(\al)}$. We assume that $t\ge t^*:=\frac 1{\sqrt 2 \al}$. Let $H$ be the halfplane containing $tD(\al)$ and having $z$ on its bounding line.

\medskip
Suppose that when the game is in position $z$, Pusher's offer is $v=(\cos \phi,\sin \phi) \in S^*$ so $|\phi|\le \frac {\pi}4$. Under the assumption that $t=t(z)\ge t^*$, Chooser's strategy gives the next position $w$ according to the following rules, see Figure~\ref{fig:Dal)}.

\smallskip
{\bf Case 1} when $z\in C_1$, her choice is $\eps=-1$ and $w=z-v.$

\medskip
{\bf Case 2} when $z\in C_2$, she chooses  $\eps \in \{-1,1\}$ so that $w=z+\eps v \in H$.

\begin{figure}[h]
\centering
\includegraphics[scale=0.8]{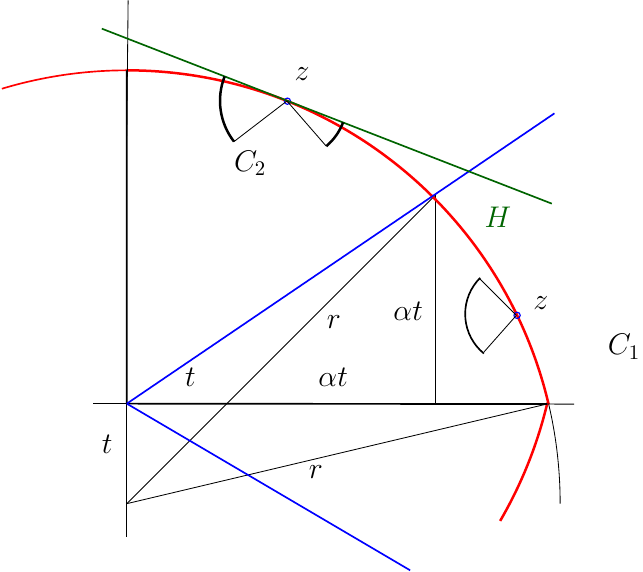}
\caption{$tD(\al)$ and the cases $z\in C_1$ and $z\in C_2$.}
\label{fig:Dal}
\end{figure}

\medskip
This strategy is not defined when $t(z) < t^*$. This is not going to matter since with $t(z)<t^*$, $z$ lies in $t^*D(\al)$ and then $\|z\|$ is bounded independently of $n$.

\medskip
We are going to estimate $t(w)=|w|_{D(\al)}$. Set $r=\sqrt 2(\al+1)t$ for $t>0$. It is clear that
\[
z \in tD(\al)=D(t,r) \subset D(t,\sqrt {r^2+1}).
\]
We want to find $s$ larger (but only slightly larger) than $t$ such that $D(t,\sqrt {r^2+1}) \subset sD(\al).$ Of course $sD(\al)=D(s,R)$ and then
$R=\sqrt 2(1+\al)s$.

\begin{figure}[h]
\centering
\includegraphics[scale=0.8]{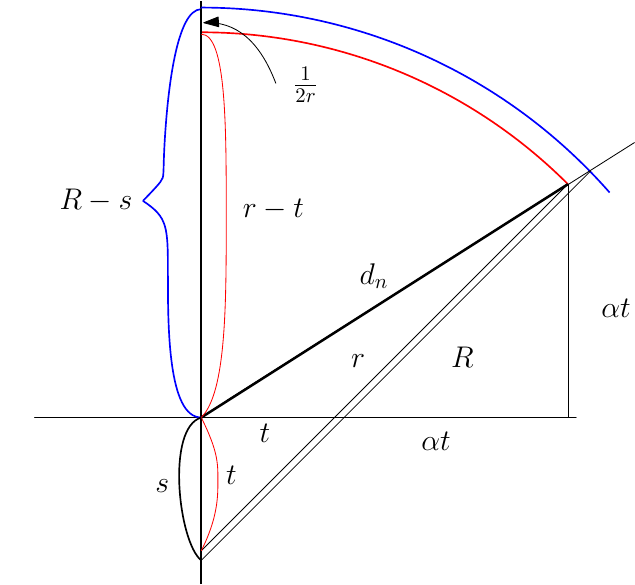}
\caption{The determination of $s$.}
\label{fig:Rands}
\end{figure}

\medskip
Observe that $\sqrt {r^2+1} < r+\frac 1{2r}$. We choose $s$ (see Figure~\ref{fig:Rands}) so that
\[
R-s =(r-t)+\frac 1{2r}
\]
.

Since by homothety $\frac {R-s}s=\frac {r-t}t$ one can check that
\begin{equation}\label{eq:t-to-s}
s =t+\frac t{2r(r-t)}=t+\frac 1{2\be t}
\end{equation}
where $\be = 2(\al+1)^2-\sqrt 2(\al+1)$. Thus $s$ is indeed only slightly larger than $t$.
The map $t \to f(t):=t+\frac 1{2\be t}$ is increasing for $t \ge \frac 1{\sqrt{2\be}}$. It is easy to check that  $t^*=\frac 1{\sqrt 2 \al}> \frac 1{\sqrt{2\be}}$.
For later reference we record here that the map
\begin{equation}\label{eq:incr} t \to f(t)=t+\frac 1{2\be t} \mbox{ is increasing for }t \ge t^*.
\end{equation}

\medskip
We show next that $t=t(z)\ge t^*$ implies that $w \in sD(\al)$, or what is the same $t(w)\le s=f(t(z)).$
\begin{claim}\label{cl:t-to-s} If $t^*\le t(z)$, then $w \in sD(\al)$ with the above choice of $s$.
\end{claim}

{\bf Proof.} By the choice of $s$, $D(t,\sqrt{r^2+1}) \subset sD(\al)$ and $z$ lies in the intersection of the upper halfplane and $B(-t,r)$.  In Case 2  $z=(x,y)$ is on the boundary of $tD(\al)$. The condition $t\ge t^*$ ensures that $y\ge t^*\al\ge \frac 1{2}$. Then $w$ also lies in the upper halfplane because its $y$ component, $y+\eps\sin \phi \ge y-\frac 1{\sqrt 2}\ge 0$. The strategy implies that $w \in B(-t,\sqrt{r^2+1}) \subset B(-t,r+\frac 1{2r})$. So $w \in sD(\al)$.


\medskip
In Case 1 we have to check that $z-v \in sD(\al)$. The critical case is when $w=z-v$ is farthest from the point $(0,-t)$, see Figure~\ref{fig:Dal}. This happens when $z=t(\al+1,\al)$ and $\phi=-\frac {\pi}4$. And even then $z-v \in B(-t,\sqrt{r^2+1})$ so $w \in sD(\al).$\qed

\medskip
Renaming $t$ as $t_k$ and $s$ as $t_{k+1}$ (\ref{eq:t-to-s}) gives rise to the following recursion (starting with $k=h$, say):
\begin{equation}\label{eq:recur}
t_{k+1}=t_k+\frac 1{2\be t_k}.
\end{equation}
In order to analyze Chooser's strategy we have understand the behaviour of this recursion.

\begin{claim}\label{cl:recur} In this recursion for every $k > h$
\[
t_h^2+\frac{k-h}{\be} \le t_k^2 <t_h^2+\frac{k-h}{\be} +\frac 1{4\be} \log \frac {k-h}h.
\]
\end{claim}

{\bf Remark.} Only the upper bound on $t_k^2$ is needed in the sequel, but the proof of Claim~\ref{cl:recur} uses the lower bound (which is actually quite close to the upper one).

\medskip
{\bf Proof.} Taking squares gives $t_{k+1}^2=t_k^2+ \frac 1{\be}+\frac 1{4\be^2 t_k^2}>t_k^2+ \frac 1{\be}$. Then
\[
t_{k+1}^2>t_k^2+\frac 1{\be}> \ldots > t_h^2+\frac {k+1-h}{\be},
\]
showing that $\frac 1{t_{k+1}^2}<\frac {\be}{k+1-h}$. So we have
\begin{eqnarray*}
t_{k+1}^2&=&t_k^2+\frac 1{\be}+\frac 1{4\be^2 t_k^2}=t_{k-1}^2+\frac 2{\be}+\frac 1{4\be^2}\left(\frac 1{t_k^2}+\frac 1{t_{k-1}^2}\right)=\ldots \\
                &=&t_h^2+\frac {k+1-h}{\be}+\frac 1{4\be^2}\left(\frac 1{t_k^2}+\ldots +\frac 1{t_h^2}\right) < \\
                &<&t_h^2+\frac {k+1-h}{\be}+\frac 1{4\be}\left(\frac 1{k-h}+\ldots + \frac 1{h+1}\right) \\
              &<&t_h^2+\frac {k+1-h}{\be} +\frac 1{4\be} \log \frac {k-h}h. \hskip3.1cm \qed
\end{eqnarray*}

\medskip
We return to Chooser's strategy now. Assume that game proceeds through positions $z_0=0,z_1,z_2,\ldots$. The set $I=\{k:t(z_k)<t^*\}$ is the union of sets of consecutive integers $\{m_i,\ldots, n_i-1\}$, $i=0,1,\ldots,$ here $m_0=0$, $n_0\ge 1$, and $m_i<n_i$. Here $n_i=\infty$ is possible but for at most one $i$. The complement of $I$, the set of subscripts $k$ with $t(z_k)\ge t^*$, is again the union of sets of consecutive integers $\{n_i,\ldots, m_i-1\}$, $i=0,1,\ldots$.  When $k \notin I$, then it is in some $\{n_i,\ldots, m_i-1\}$. Since $t(z_{n_i-1}) < t^*$, $z_{n_i}\in z_{n_i-1}+B \subset t^*D(\al)+B$. Set $t_*=\max \{t(u): u\in t^*D(\al)+B\}$, a constant only depending on $\al$. Moreover $t(z_{n_i})\le t_*$ for every $i$.

\begin{claim}\label{cl:tk} With this strategy $t(z_k)^2<t_*^2+\frac k{\be}+\frac {\log k}{4\be}$ for every $k\ge 1$.
\end{claim}

{\bf Proof.}  For $k \in I$, $t(z_k)<t^* < t_*$ so we are left with the case $k \notin I$. Then $k \in J:=\{n_i,\ldots,m_i-1\}$ for some $i$.

\medskip
The recursion~(\ref{eq:recur}) starting at $h=n_i$ with $t(z_h)=t_h\ge t^*$ gives numbers $t_k$ for every $k \in J,$ that depend on $h$ but that matters not. The upper bound in the previous claim shows that for $k \in J$ (except $k=h$)
\[
t_k^2 < t_h^2+ \frac {k-h}{\be}+\frac 1{4\be}\log \frac {k-h}{h} <t_*^2 + \frac k{\be}+\frac {\log k}{4\be}
\]
because $t(z_h)=t_h\le t_*$ and $h\ge 1$.

\medskip
Using (\ref{eq:incr}), an easy induction on $k$ starting with $k=h$ shows that $t(z_k)\le t_k$ for every for $k \in J$.
Consequently for every $k>0$, $t(z_k)^2 < t_*^2 +\frac {k}{\be}+\frac {\log k}{4\be}$.\qed

\bigskip
\section{Proof of the upper bound in Theorem~\ref{th:main}}\label{sec:proof}

Define $T_n>0$ via $T_n^2= t_*^2 +\frac {n}{\be}+\frac {\log n}{4\be}$. 
Write $d_n=\|T_n(1+\al,\al)\|$. Lemma~\ref{cl:tk} gives that
\[
d_n < \sqrt n \sqrt{\frac{\al^2+{(\al+1)^2}}{\be}}\left(1+O\left(\frac {\log n}{\sqrt n}\right)\right).
\]
We are to choose $\al$ so that
$$\frac{\al^2+{(\al+1)^2}}{\be}=\frac{\al^2+{(\al+1)^2}}{2(\al+1)^2-\sqrt 2(\al+1)}$$
is minimal. The usual computation shows that this happens when $\al=\frac 1{\sqrt 2}+\sqrt{1+\frac 1{\sqrt 2}}=2.013669..$ and with this $\al$,  $\be=13.895312..$ and $d_n$ is approximately $0.972112.. \sqrt n$.

\medskip
We fix $\al=\frac 1{\sqrt 2}+\sqrt{1+\frac 1{\sqrt 2}}=2.013669..$ and then $t^*=0.351272..$. Set $D=D(\al)$ with this $\al$,  see Figure~\ref{fig:DB}.
A simple computation shows that $t_*=0.511187..$.

\medskip
The proof of the upper bound in Theorem~\ref{th:main} is finished by the following lemma.

\begin{lemma}\label{l:fin} $\|z_n\| <d_n+1.$
\end{lemma}

{\bf Proof.} We prove the stronger statement that
\[
\max_{j\in \{1,\ldots, n\}} \|z_j\| <d_n+1.
\]
Assume this maximum is reached on $z_k=(x,y)$ and $x,y\ge 0$ by symmetry. By Claim~\ref{cl:tk} $t(z_k) < T_k.$

\medskip
If $z_k \in C_2$,  then $z_k\in T_kD\cap C_2$ whose farthest point from the origin is $T_k(1+\al,\al)$ (see Figure~\ref{fig:Rands}) so $\|z_k\|\le d_k\le d_n$.

\medskip
So $z_k \in C_1$ and assume that $z_{k-1} \in C_1$ too. Then by Chooser's strategy $z_k=z_{k-1}-v_{k-1}$ where $v_{k-1}=(\cos \phi, \sin \phi)$ with $|\phi|\le \frac {\pi}4$, and $\|z_{k-1}\| > \|z_k\|$ contradicting the maximality of $\|z_k\|$. So $z_{k-1} \in C_2$ and $z_{k-1} \in T_{k-1}D$, as we just checked,  $\|z_{k-1}\|\le d_{k-1}.$ Then $\|z_k\| \le \|z_{k-1}\|+1\le d_{k-1}+1<d_n+1$.\qed

\section{Prologue}\label{sec:pro}

\medskip
After hearing the above results Joel Spencer asked me what happens in game $(S^{\ga},B)$ where $S^{\ga}=
\{(\cos \phi,\sin \phi): |\phi|\le \ga\}.$ The previous methods work in this case up to a certain level as described below.

\medskip
For the lower bound set $K^{\ga}=D(\cot \ga, (\sin \ga)^{-1})$ and check that in game $(S^{\ga},B)$,
$|v_1|_{K^{\ga}}\ge 1$, just like in (\ref{eq:v_1}). The strategy of Pusher is the same as in game $(S^*,B)$:
in step $n$ he computes $t=|z_n|_{K^{\ga}}$ and selects a unit vector $v_n\in S^{\ga}$ tangent to $tK^{\ga}$ at $z_n\in \bd (t_nK^{\ga})$.
(Again such a unit vector always exists.) Assume that the $y$ component of $z_n$ is non-negative. Then $z_n$ is at distance $r=t(\sin \ga)^{-1}$ from $(0,-t\cot \ga)$. Now $z_{n+1}=z_n+\eps_nv_n$, and for both $\eps= \pm 1$, $z_{n+1}$ is exactly at distance $\sqrt{r^2+1}$ from $(0,-t\cot \ga)$. Direct checking shows again that $z_{n+1}$ does not lie in the interior of $\sqrt{t^2+1}K^{\ga}$, implying that $t(z_{n+1})\ge \sqrt {t(z_n)^2+1}$. This leads to the recursion $t(z_{n+1})^2\ge t(z_n)^2+1$ which shows again that $t(z_n)^2\ge n.$ This strategy gives that
\begin{equation}\label{eq:upper}
\|z_n\|\ge (1-\cos \ga)\sqrt n= \frac {\sin^2 \ga}{1+\cos \ga}\sqrt n.
\end{equation}
So for small $\ga>0$ Pusher can guarantee that $\|z_n\|$ is of order $\sqrt n$ but only with a small multiplicative factor.

\medskip
The previous strategy for Chooser works again but only when $\ga$ is below $\ga_0= 0.7967..$ in radian (or $45.64^{\circ}$) as we shall see soon. We assume from now on that $\ga <0.85$.

\begin{figure}[h]
\centering
\includegraphics[scale=0.8]{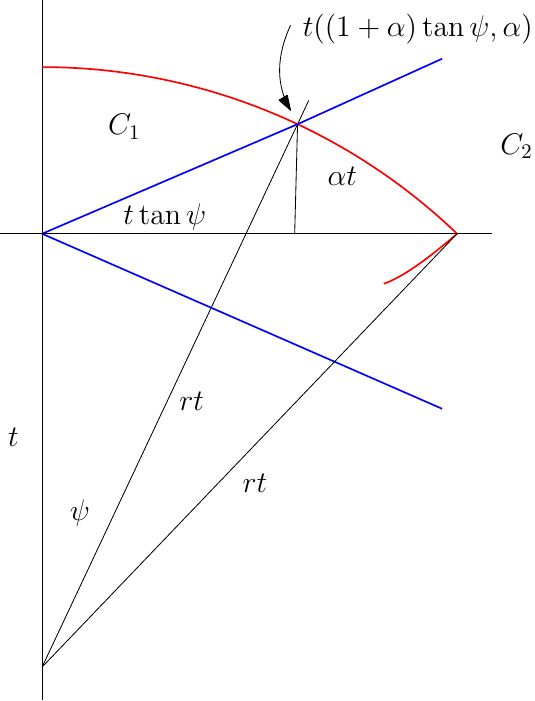}
\caption{The top right quarter of $tD^{\ga}$}
\label{fig:Dpsi}
\end{figure}

This time we define $D^{\ga}:=D^{\ga}(\al)=D(1,r)$ where $r=\frac {(1+\al)}{\sin \ga}$ with $\al> 0$ that will be fixed later.
The cones $C_1$ and $C_2$ are defined analogously, see Figure~\ref{fig:Dpsi} where the angle $\psi$ is determined by $\al$ and  $\ga$.

\medskip
Suppose $z=(x,y)$ is the current position in game $(S^*,D^{\ga})$. We assume that $x,y\ge 0$. Set $t=t(z)\!:=|z|_{D^{\ga}}$. We assume that $t\ge t^*=\frac {\cos \ga}{\al}$. Let $H$ be the halfplane containing $tD^{\ga}$ and having $z$ on its bounding line.

\medskip
Assume that Pusher's offer is $v=(\cos \phi,\sin \phi)$ with $|\phi|\le \ga$. Chooser uses the same strategy to determine the next position $w$ as before:

\smallskip
{\bf Case 1} when $z\in C_1$, her choice is $\eps=-1$ and $w=z-v.$

\smallskip
{\bf Case 2} when $z\in C_2$, she chooses  $\eps \in \{-1,1\}$ so that $w=z+\eps v \in H$.

\medskip
This strategy is not defined for $t < t^*$ but that is no problem since when $t(z)<t^*$, $z$ lies in $t^*D^{\ga}$ and $\|z\|$ is bounded by a function of $\al$ and $\ga.$ The analysis is the same as before: set
\[
s=t+\frac t{2rt(rt-t)}=t+\frac 1{2\be t} \mbox{ where } \be =\frac {1+\al}{\sin \ga}\left(\frac {1+\al}{\sin \ga}-1\right).
\]
The map $t \to f(t)=t+\frac 1{2\be t}$ is increasing for $t>\frac 1{2\be}$ so it is increasing for $t\ge t^*$ because $t^*>\frac 1{2\be}$.
The analogue of Claim~\ref{cl:t-to-s}, that is, $w \in sD^{\ga}$ with this choice of $s$ is proved the same way. In this case the recursion (\ref{eq:recur}) takes an analogous form and Claim~\ref{cl:tk} holds again (details omitted):
\[
t_k^2<t_*^2+\frac k{\be}+\frac 1{4\be}\log k \mbox{ for every }k\ge 0.
\]

Returning to the sequence of positions $z_0=0,z_1,\ldots$ we have again $t(z_n)\le T_n$.
Writing $d_n$ for the distance between the origin and $(t(1+\al)\tan \ga,\al)$ we see that
\begin{eqnarray*}
d_n^2&\le \frac n{\be}\left((1+\al)^2\tan^2 \ga +\al^2\right)\left(1+O(\log n)\right)\\
         &=n \sin^2\ga\,  G(\al,\ga)\left(1+O(\log n)\right).
\end{eqnarray*}
where the function $G(\al,\ga)$ is defined for all $\al>0$ and $0.85 \ge \ga>0$ by
\[
G(\al,\ga)=\frac {(1+\al)^2\tan^2 \ga +\al^2}{(1+\al)(1+\al-\sin \ga)}.
\]

For a fixed $\ga \in (0,0.85)$ we want to see if there is an $\al>0$ such that the function $F(\al,\ga):=\sin^2\ga\,  G(\al,\ga)<1$ as otherwise $d_n\ge \sqrt n(1+o(1))$ and the simple upper bound $\sqrt n$ from game $(S^1,B)$ is either better or is just as good. A calculation (using Mathematica) shows that when $\ga > \ga_0$, $F(\al,\ga)$ is larger than 1 for all $\al>0$, where $\ga_0 = 0.7967..$, so for that range this method does not give an upper bound better than $\sqrt n.$ It is not clear to me if there is a strategy for Chooser that would give $\|z_n\|<c\sqrt n$ for some $c<1$ in the range $\ga > \ga_0$.

\medskip
On the other hand, for fixed $\ga<\ga_0$, $\min_{\al>0}F(\al,\ga)<1$ so the above strategy for Chooser gives an upper bound
on $\|z_n\|$ better than $\sqrt n$. We made no effort to determine the minimizer $\al$ in this case except when $\ga$ is close to $0$. Choosing for instance $\al=\ga^{1.5}$ gives the bound
\[
\|z_n\|\le \sin^2 \ga \sqrt n(1+o(1)),
\]
which is not far from the lower bound in (\ref{eq:upper}).

\bigskip
{\bf Acknowledgements.} This piece of work was partially supported by NKFIH grants No 131529, 132696, and 133919 and also by the HUN-REN Research Network. Part of this research was carried out on a very pleasant and fruitful visit to the Hausdorff Institute of Mathematics in Bonn whose support is also acknowledged. I am indebted to Gergely Ambrus for useful discussions at the starting stages of this research.

\vskip0.3cm

\noindent
Imre B\'ar\'any \\
Alfr\'ed R\'enyi Institute of Mathematics, HUN-REN\\
13 Re\'altanoda Street, Budapest 1053 Hungary\\
{\tt barany.imre@renyi.hu}\\
and \\
Department of Mathematics\\
University College London\\
Gower Street, London, WC1E 6BT, UK\\

\end{document}